\documentclass{elsarticle}

\usepackage{setspace}
\usepackage{amsmath}
\usepackage{lipsum}
\usepackage{bm}
\usepackage{amsfonts}
\usepackage{graphicx}
\usepackage{epstopdf}
\usepackage[ruled,vlined]{algorithm2e}
\usepackage{algorithmic}
\usepackage{dirtytalk}
\usepackage{graphicx,epstopdf}
\usepackage[caption=false]{subfig}
\usepackage{tikz}
\usepackage{pgfplots}
\pgfplotsset{compat = 1.3}

\journal{Journal of Computational Physics}

\begin{document}

\begin{frontmatter}

\title{Agglomeration-Based Geometric Multigrid Solvers for Compact Discontinuous Galerkin Discretizations on Unstructured Meshes}

\author[rvt2,rvt3]{Y.~Pan\fnref{fn1}\corref{cor1}}
\ead{yllpan@berkeley.edu}

\author[rvt2,rvt3]{P.-O.~Persson\fnref{fn3}}
\ead{persson@berkeley.edu}

\address[rvt2]{Department of Mathematics, University of California, Berkeley, Berkeley, CA 94720, United States}
\address[rvt3]{Mathematics Group, Lawrence Berkeley National Laboratory, 1 Cyclotron Road, Berkeley, CA 94720, United States}
\cortext[cor1]{Corresponding author}
\fntext[fn1]{Graduate student, Department of Mathematics, University of California, Berkeley}
\fntext[fn3]{Professor, Department of Mathematics, University of California, Berkeley}

\begin{keyword}
 discontinuous Galerkin, %
 agglomeration, %
 geometric multigrid
\end{keyword}

\begin{abstract}
  We present a geometric multigrid solver for the Compact Discontinuous Galerkin method through building a hierarchy of coarser meshes using a simple agglomeration method which handles arbitrary element shapes and dimensions. The method is easily extendable to other discontinuous Galerkin discretizations, including the Local DG method and the Interior Penalty method. We demonstrate excellent solver performance for Poisson's equation, provided a flux formulation is used for the operator coarsening and a suitable switch function chosen for the numerical fluxes.
\end{abstract}

\end{frontmatter}

\section{Introduction}

The discontinuous Galerkin (DG) method with high-order approximations are becoming increasingly popular for the solution of systems of conservation laws, due to their natural ability to stabilize convection-dominated problem on arbitrary unstructured meshes with high-order accuracy. The resulting semi-discrete systems are often integrated in time using explicit solvers, however, for many real-world problems it is widely believed that implicit solvers will be required. This poses many challenges, since the Jacobian matrices are expensive to compute and store, and specialized solvers are required to solve the corresponding linear systems that arise.

One of the most important solver techniques employed, at least for elliptic or diffusion-dominated problems, is the multigrid method \cite{brandt77multilevel}. The method has been used extensively for DG methods \cite{fidkowski05multigrid,persson08newtongmres,biros15multigrid,pazner20precondition,saye19multigrid,saye19interface}, where it can naturally be applied as a $p$-multigrid solver where the grid hierarchy is formed by varying the polynomial degrees in each element. It can also be used in the more traditional $h$-multigrid setting, where the hierarchy is based on meshes of varying coarseness, or as a combined $hp$-multigrid method which combines both these techniques \cite{mavriplis07paralleldg}. For fully unstructured meshes, it is in general difficult to coarsen a given mesh in order to produce the mesh hierarchies needed for a full $h$-multigrid. This is one of the motivations for using so-called Algebraic Multigrid methods \cite{xu17amg,boomerAMG,mfem}. 

An alternative approach for coarsening an unstructured mesh is agglomeration, that is, merging neighboring elements into larger ones successively. The technique is not widely used for continuous Galerkin finite element methods, because of the difficulties in defining continuous approximation spaces on the resulting polyhedral elements. However, the technique has been used successfully for finite volume methods \cite{vassilevski01agglomeration,chan98agglomeration,strauss03agglomeration}, where it is easier to update the element-averages after coarsening. This is also true for high-order discontinuous Galerkin methods, since they are straight-forward to implement on meshes of arbitrarily shaped elements \cite{berggren10agglomeration,dargaville20agglomeration}.

In this paper, we propose an agglomeration-based $h$-multigrid method for Poisson's equation based on the CDG method \cite{peraire08cdg}. This is a variant of the LDG method \cite{cockburn98ldg}, with important benefits such as element-wise compact stencils and improved stability properties. However, our method should be straight-forward to use with the LDG method, or any other discretization such as Interior Penalty or the BR2 methods \cite{arnold02unified}.

We perform the element agglomeration with a simple approach which extends to arbitrary elements and dimensions. While the resulting hierarchy might not be optimal for the multigrid performance, our numerical experiments demonstrate that the method is quite insensitive to the shape of the agglomerated elements. We also show the importance of choosing a good switch function for the numerical fluxes in the CDG method.

The paper is organized as follows. In Section 2, we describe the CDG discretization and in particular write it in the so-called flux formulation which is needed for the operator coarsening in the multigrid method. Next, we outline the (heuristic) geometric element agglomeration algorithm in Section 3, and the details of the multigrid method in Section 4. Our numerical results in Section 5 show a number of important properties of our scheme, and demonstrate its performance.

\section{Discontinuous Galerkin formulation}
\label{sec:main}
\subsection{Problem definition}
For our notation, quantities that have a spatial dimension, such as the spatial gradient of a function, are bolded whilst scalar functions are not. We consider here Poisson's equation as our model elliptic problem 
\begin{equation} \label{eq:poisson}
    \begin{aligned}
        \nabla^2 u &= f &\text{in } &\Omega, \\
        u &= g_D &\text{on } &\Gamma_D, \\
        \nabla u \cdot \bm{n} &= g_N &\text{on } &\Gamma_N,
    \end{aligned}
\end{equation}
in a domain $\Omega \subset \mathbb{R}^d$, where $d \in \{ 1,2,3 \}$ is the dimension of the system. $\Gamma_D, \Gamma_N$ respectively denote parts of the boundary $\partial \Omega$ on which Dirichlet and Neumann boundary conditions are imposed, with $\bm{n}$ denoting the unit outward normal on $\partial \Omega$. Here, $f(\bm{x})$ is an arbitrary given function in $L^2(\Omega)$ and we further assume that the length of $\Omega_D$ is strictly greater than zero.  

\subsection{DG formulation for elliptic problems}
To apply a DG method to the above model problem, we rewrite Equation \ref{eq:poisson} as a first order system of equations by introducing the variable $\bm{q} = \nabla u$ and rewriting the Laplacian operator as the divergence of $\bm{q}$,

\begin{equation} \label{eq:poisson_ldg}
    \begin{aligned}
        \nabla \cdot \bm{q} &= f  &\text{in } &\Omega, \\
        \bm{q} &= \nabla u &\text{in } &\Omega, \\
        u &= g_D &\text{on } &\Gamma_D, \\
        \bm{q} \cdot \bm{n} &= g_N &\text{on } &\Gamma_N.
    \end{aligned}
\end{equation}

In this work, we consider discretizations where meshes $\mathcal{T}_h = \{ K \}$ of $\Omega$ may consist of arbitrarily shaped elements, with the only restriction being elements must not self intersect. We define the broken spaces $V(\mathcal{T}_h)$ and $\Sigma(\mathcal{T}_h)$ as the union of Sobolev spaces $H^1(K)$ and $[H^1(K)]^d$ restricted to each element $K$. Specifically,

\begin{align}
    V &= \{ v \in L^2(\Omega) : v|_{K} \in H^1(K), ~\forall K \in \mathcal{T}_h \}\\
    \Sigma &= \{ \bm{\tau} \in [L^2(\Omega)]^d : \bm{\tau}|_{K} \in [H^1(K)]^d, ~\forall K \in \mathcal{T}_h \}
\end{align}

We also introduce the finite element spaces $V_h \subset V$ and $\Sigma_h \subset \Sigma$ as
\begin{align}
    V_h &= \{ v \in L^2(\Omega) : v|_{K} \in \mathcal{P}_p(K), ~\forall K \in \mathcal{T}_h \}\\
    \Sigma_h &= \{ \bm{\tau} \in [L^2(\Omega)]^d : \bm{\tau}|_{K} \in [\mathcal{P}_p(K)]^d, ~\forall K \in \mathcal{T}_h \}
\end{align}
where $\mathcal{P}_p(K)$ denotes the space of polynomial functions of order at most $p \geq 1$ on each element $K$.

We obtain a weak DG formulation by multiplying the system of equations with test functions $v, \bm{\tau}$ before integrating by parts. From this our formulation can be expressed as finding $\bm{u}_h \in V_h, \bm{q}_h \in \Sigma_h$ such that for all $K \in \mathcal{T}_h = \{ K \}$, we have

\begin{equation} \label{eq:poisson_weak}
    \begin{aligned}
        \int_K (\bm{q}_h + u_h \nabla) \cdot \bm{\tau} dx = \int_{\partial K} \hat{u} \bm{\tau} \cdot \bm{n} ds \quad &\forall \bm{\tau} \in [\mathcal{P}_p(K)]^d, \\
        \int_K \bm{q_h} \cdot \nabla v ~ dx = \int_{\partial K} v \hat{\bm{q}} \cdot \bm{n} ds + \int_{K} fv ~ dx \quad &\forall v \in \mathcal{P}_p(K).
    \end{aligned}
\end{equation}

The numerical fluxes $\hat{u}, \hat{\bm{q}}$ approximate the quantities to $u$ and to $\bm{q} = \nabla u$ on the boundaries of each element $K$. For the CDG method, numerical fluxes are expressed as a function of the fields $u_h$ and $\bm{q}_h$, in addition to the specified boundary conditions on $\partial \Omega$ as follows.

To specify the numerical fluxes, we define a switch function $S_{K}^{K'} \in \{ -1,1 \}$ on each internal boundary separating element $K$ from its neighbour $K'$, which satisfies the property $S_{K}^{K'} = -S_{K'}^{K}$. One example is the natural switch function, where given any enumeration of the elements $\{ \mathcal{N}(K) \}$, for any two elements $K,K'$, the switch $S_{K}^{K'} > 0$ if $\mathcal{N}(K) > \mathcal{N}(K')$. Given a switch function, the numerical fluxes are defined as:

\begin{itemize}
    \item In Equation \ref{eq:poisson_weak}, $\hat{u}$ is defined by standard upwinding based on the switch function
    \begin{equation}
        \hat{u} = 
        \begin{cases}
          u_h &\text{if } S_{K}^{K'} > 0 \\
          u'_h     &\text{if } S_{K}^{K'} < 0
        \end{cases}
    \end{equation}
    where $u'_h$ is the numerical solution to $u$ in Equation \ref{eq:poisson_weak} on the neighbouring element $K'$ on boundary $\partial K$.
    \item On every inter-element boundary $f$ separating two elements $K, \tilde{K}$, where $S_{K}^{\tilde{K}} < 0$, define a \say{boundary gradient} $\bm{q}^f_h$ using a slight modification of Equation \ref{eq:poisson_weak}
    \begin{equation} \label{eq:edge_flux_interior}
        \int_{K} (\bm{q}_h^f + u_h \nabla) \cdot \bm{\tau} dx = \int_{\partial K \setminus f} u_h \bm{\tau} \cdot \bm{n} ds + \int_{f} \tilde{u}_h \bm{\tau} \cdot \bm{n} ds
    \end{equation}
    where tilde on $\tilde{u}_h, \tilde{\bm{q}}_h$ denotes numerical solutions to the respective fields defined on $\tilde{K}$. The flux $\hat{\bm{q}}$ on $f$ is then defined simply by restricting $\bm{q}_h^f$ to the boundary $f$.
    \item On a boundary $f$ of element $K$ that coincides with $\partial \Omega$, we similarly define a \say{boundary gradient}
    \begin{equation} \label{eq:edge_flux_boundary}
        \int_{K} (\bm{q}_h^f + u_h \nabla) \cdot \bm{\tau} dx = \int_{\partial K} u \bm{\tau} \cdot \bm{n} ds
    \end{equation}
    The numerical fluxes are defined using the defined \say{boundary gradient} in addition to the specified boundary conditions,
    \begin{equation}
        \label{eq:cdg_fluxes}
        \begin{aligned}
            \hat{\bm{q}} &= \bm{q}_h^f - C_{D}(u_h - g_D)\bm{n}, &\hat{u} = g_D &~\text{ on } \partial \Omega_D \\
            \hat{\bm{q}} &= g_N \bm{n}, &\hat{u} = u_h &~\text{ on } \partial \Omega_N
        \end{aligned}
    \end{equation}
    where the parameter $C_D > 0$ is included for additional stabilisation. For our applications, we choose $C_D = \gamma/h_{avg}$, where $\gamma > 0$ is a constant, and $h_{avg}$ is mean height of elements $K$ on the boundary $\partial \Omega_D$. The choice of $\gamma$ and its effect on multigrid convergence is discussed in Section \ref{sect:fluxprimal}.
\end{itemize}
We briefly note the similarity of the CDG method to the LDG method, with the only distinction being in the definition of fluxes $\hat{\bm{q}}$. For a more detailed treatment on the CDG method and its properties we turn the reader to \cite{peraire08cdg}.

\subsection{Discrete formulation}
Discretising Equation \ref{eq:poisson_weak} we obtain a linear system
\begin{equation}
    \begin{split}
        M\bm{q}_h + Gu_h &= \bm{r} \\
        \tilde{D}\bm{q}_h + N\hat{\bm{q}} &= f
    \end{split}
\end{equation}
where $M$ denotes the system mass matrix, $G$ the discrete gradient operator, and $\bm{r}$ the Dirichlet vector, defined as
\begin{equation}
    \label{eq:disc_grad}
    \begin{split}
         G(u_h) &= \sum_K \Bigg( \int_K u_h \nabla \cdot \bm{\tau} dx - \int_{\partial K \setminus \partial \Omega_D} \hat{u} \bm{\tau} \cdot \bm{n} ds \Bigg) \\
         \bm{r} &= \int_{\partial \Omega_D} g_D \bm{\tau} \cdot \bm{n} ds
    \end{split}
\end{equation}
The operators $\tilde{D}$, $N$, and the vector $f$ are defined as
\begin{equation}
    \begin{split}
        \tilde{D}(\bm{q}_h) &= \sum_K \int_{K} \bm{q}_h \cdot \nabla v dx \\
        N(\hat{\bm{q}}) &= -\sum_K \int_{\partial K \setminus \partial \Omega_N} v \hat{\bm{q}} \cdot \bm{n} ds \\
        f &= \sum_K \int_K fv dx - \int_{\partial \Omega_N} v g_N ds - \int_{\partial \Omega_D} v C_D g_D ds
    \end{split}
\end{equation}

Following \cite{peraire08cdg}, it is possible to write $\hat{\bm{q}}$ in terms of the variables $u_h$ and $\bm{q}_h$ as a consequence of Equation \ref{eq:cdg_fluxes}. Similarly, following Equation \ref{eq:disc_grad}, it is possible to write the variable $\bm{q}_h$ simply as a variable of the unknown $u_h$. All together this allows us to write the Poisson system discretely as a single linear system
\begin{equation} \label{eq:primal_form}
    Au_h = b
\end{equation}
where the matrix $A$ is compact, meaning a block $(i,j)$ has non-zero entries if and only if elements $K_i, K_j$ are immediate neighbours. This is known as the primal form of the CDG method.

For construction of a multigrid solver however, following \cite{saye19multigrid}, direct coarsening of the primal operator $A$ can lead to decreased performance of the solver. A flux form of the CDG system can instead be defined as follows
\begin{equation} \label{eq:flux_matrix_form}
    \begin{bmatrix} M & G \\ D & C \end{bmatrix}
    \begin{bmatrix} \bm{q}_h \\ u_h \end{bmatrix} = 
    \begin{bmatrix} \bm{r} \\ s \end{bmatrix}
\end{equation}
where $D=-G^T$ is the discrete divergence operator, and $A = C - DM^{-1}G, ~s = f - DM^{-1}\bm{r}$. Unlike the matrix $A$, the matrix $C$ is not compact in that it may contain non-zero entries in a block $(i,j)$ where elements $K_i, K_j$ are not immediate neighbours. 

The CDG method is closely related to the LDG method in that the only difference in the flux formulations of the two lie in the bottom right hand entry of the flux operator; for LDG the $C$ matrix is equal to the zero matrix. The $C$ matrix in the CDG method cancels out the non-compact entries from the term $DM^{-1}G$, rendering the resulting matrix $A$ to be compact. This therefore implies that in general it is unnecessary to store the non-compact entries of the matrix $C$ as they may be implicitly inferred.

\section{Mesh hierarchy}
For h-multilevel solvers, a hierarchy of successively coarser mesh discretizations are constructed on which the matrix equation $Ax=b$ is solved approximately on each level of the hierarchy. While for structured meshes coarsening algorithms such as quadtree/octree are widely used to construct mesh hierarchies, it is in general difficult to construct such hierarchies for unstructured meshes.

\subsection{Element agglomeration}
In this work, we focus on mesh hierarchy construction via agglomeration. A valid mesh agglomerate is defined in this work as a partition of the set of elements $\mathcal{T}_h = \{ K \}$ such that the union of elements within each subset of the partition form a connected domain. This choice of hierarchy constructed is based on the observation that in general for a DG formulation, the lack of $C^0$ continuity required in standard continuous finite elements allows for the easy definition of modal basis functions on arbitrarily shaped polyhedra.

\begin{figure}
\centering 
    \label{fig:meshpartions}\includegraphics[scale=0.4]{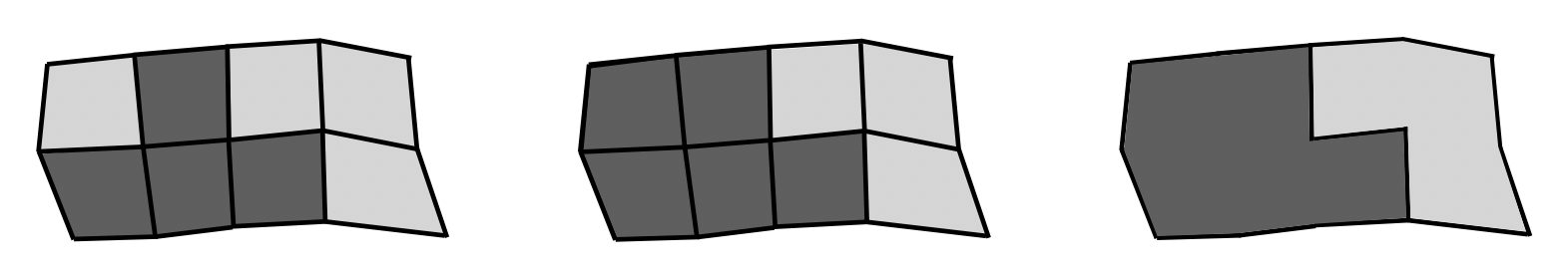}
    \caption{Example mesh partitions. The left partition is invalid as one subset is not connected. The middle partition is valid and agglomerated into the polygonal mesh on the right.}
\end{figure}

To define a mesh hierarchy, elements within each partition are agglomerated to form a single polygonal element, which are then all collected as the set of elements for the next level of the mesh hierarchy. This process can be performed recursively until the final level of the hierarchy contains only a single polygonal element defined by the boundaries of the computational domain $\Omega$.

\begin{algorithm}[H] \label{alg:vcycle}
\caption{Recursive mesh agglomeration}
    \begin{algorithmic}
        \STATE{Input mesh with elements $E_0$}
        \STATE{Hierarchy storage $E = \{ E_0 \}$}
        \STATE{$E_{level}= E_0$}
        \WHILE{$length(E_{level}) > 1$}
            \STATE{ $E_{level+1}$ = Find set partition of $E_{level}$ }
            \STATE{ Append $E_{level+1}$ to $E$ }
            \STATE{ $E_{level} = E_{level+1}$ }
        \ENDWHILE
    \RETURN $E$
    \end{algorithmic}
\end{algorithm}

The problem of finding mesh partitions is well studied in literature, including popular domain decomposition methods in the software package METIS \cite{Karypis99metis}. In this work however we use a simple greedy heuristic to demonstrate the generality of the method for mesh partitions of arbitrary shape and quality.

\subsection{Greedy agglomeration} 
\label{sect:greedy_agglomeration}
\begin{figure}
\centering 
    \includegraphics[scale=0.3]{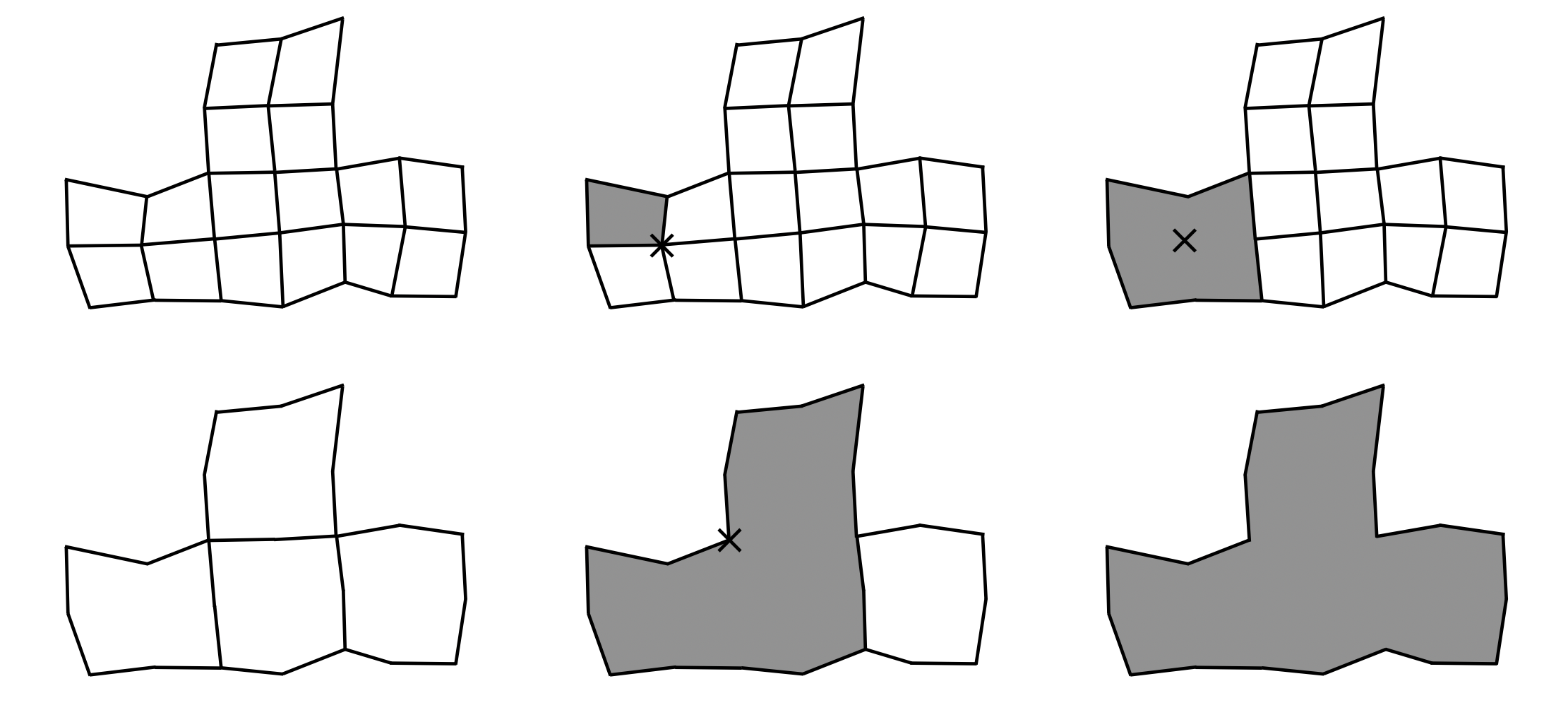}
    \caption{Schematic of agglomeration algorithm. Top row from left to right: Input mesh at level 0 of hierarchy, element with fewest number of neighbours and its vertex touching most unprocessed neighbours are chosen, all unprocessed elements touching chosen vertex are grouped to form 8 sided polygonal element for next level in hierarchy. Bottom row from left to right: Mesh at level 1 or hierarchy, all elements touching vertex are merged to form 16 sided polygonal element leaving the unshaded element with 0 unprocessed neighbours, unprocessed element with 0 unprocessed neighbour elements is merged into neighbouring agglomerate polygon. }
    \label{fig:agglomeration}
\end{figure}

We describe a mesh agglomeration algorithm through use of a simple greedy heuristic, outlined in Algorithm \ref{alg:agglomeration}. To construct a new mesh at a lower level, we assign to each element of the input mesh an integer weight corresponding to the number of neighbour elements in the mesh not yet processed. Elements are loaded into a priority queue and processed in ascending order according to the integer weights. To process an element, we identify the vertex of the element adjacent to the most unprocessed elements left in the priority queue, breaking ties at random. All the unprocessed elements touching the identified vertex are marked as processed, and agglomerated into a subset of the mesh partition, the union of which serves as an polyhedral element in the new mesh, termed a block. The priority queue is updated to reflect the removal the corresponding elements, and the algorithm repeated until no elements are left remaining in the priority queue. This process is shown in the first row of Figure \ref{fig:agglomeration}.

In the case where an unprocessed element neighbours fewer than two unprocessed elements left in the priority queue, we instead append the element to the smallest adjacent block, breaking ties at random. The priority queue is then updated to reflect the successful processing of the element. This special case is shown in the bottom row of Figure \ref{fig:agglomeration}.

Due to the use of a priority queue, which is inserted into and the minimum extracted from $n$ times respectively, where $n$ denotes the number of vertices in the mesh, the overall computational cost of the algorithm scales as $O(n\log n)$. The memory cost of the algorithm however scales only as $O(n)$, as only a single integer indicating the number of unprocessed neighbours to each vertex is stored in the priority queue.

\begin{algorithm} \label{alg:agglomeration}
\caption{Greedy mesh agglomeration}
    \begin{algorithmic}
        \STATE{Input mesh $M$ with vertices $\{v_i\}$ and elements $\{e_i\}$ }
        \STATE{Create priority queue $P = \{ (e_i, N(e_i)), \forall e_i \in M \}$, $N(e_i)=\#$neighbours of $e_i$}
        \STATE{Create empty array $N$ to store new elements}
        \WHILE{$length(P) > 0$}
            \STATE{Pop $e_{min}$ with smallest $N(e_{min})$ from $P$}
            \IF{$N(e_{min}) \geq 2$}
                \STATE{Find vertex $v_{max}$ of $e_{min}$ adjacent to the most elements in $P$ }
                \STATE{Create set $E = \{ e_j \in P, e_j \text{ adjacent to } v_{max} \} \cup \{ e_{min} \}$}
                \FOR{ $e_j$ in $E$ }
                    \STATE{Remove $e_j$ from $P$}
                    \STATE{Update all neighbours $e_k$ of $e_j$ in $P$, $N(e_k) = N(e_k)-1$}
                \ENDFOR
                \STATE{Combine all elements in $E$ to form new element, append to $N$}
            \ELSE
                \STATE{Create $N_{adj} = \{$ elements $n_j \in N$ with subelement $e_k$ adjacent to $e_{min} \}$}
                \STATE{Find $n_{min}$, element in $N$ with fewest subelements $e_k$}
                \STATE{Append $e_{min}$ to $n_{min}$}
                \STATE{Update all neighbours $e_k$ of $e_{min}$ in $P$, $N(e_k) = N(e_k)-1$}
            \ENDIF
        \ENDWHILE
    \RETURN $N$
    \end{algorithmic}
\end{algorithm}

\subsection{Basis functions and quadrature}
To define basis functions on the generated polyhedral blocks at each level in the hierarchy, we adopt a modal basis set of polynomials on each block due to the difficulty of assigning nodal basis functions on arbitrary polyhedra. For instance the linear set of basis functions of this form would be simply $1,x,y$ in two dimensions. To numerically integrate on each of the blocks, we use the fact that each of the blocks are constructed by taking a union of a subset of elements $\{ K_{n_1},...,K_{n_j} \}$ from the input mesh. This allows quadrature on polygonal elements to be computed by summing contributions from each sub-element of the block, which can be calculated using preexisting quadrature defined on the input mesh. Thus no additional computational expense due to quadrature is required at each coarser level of the mesh hierarchy.

\subsection{Solution transfer} \label{sect:sol_transfer}
We define restriction and interpolation operators to transfer residuals and states between neighbouring levels in the mesh hierarchy. For the purposes of our preconditioning strategy, we focus only on residual restriction and the case of state prolongation but do not consider the case of state restriction.

The prolongation operator from level $l+1$ to level $l$ acts as
\begin{equation}
    L_{l+1}^l v_h^{l+1} = v_h^{l}.
\end{equation}
As the basis functions for each polygonal block are chosen to be the same modal polynomials at each level, the operator can be chosen to be simple injection \cite{brenner}. Following \cite{biros10} this has an equivalent variational formulation, which can be defined using an $L^2$ projection.

The restriction operator is defined as the adjoint of the prolongation operator $L^l_{l+1}$
\begin{equation}
    (R_l^{l+1} u_h^{l}, v_h^{l+1})_{l+1} = (u_h^{l}, L_{l+1}^l v_h^{l+1})_{l}
\end{equation}
for all $u_h^{l}, v_h^{l+1}$ piecewise polynomial functions defined on levels $l, l+1$ respectively. Equivalently using the $L^2$ weak formulation the restriction operator can be written as
\begin{equation}
    M^{l+1} R_l^{l+1} = ( M^{l} L_{l+1}^l )^T
\end{equation}
where $M^{l}, M^{l+1}$ denote mass matrices for the corresponding superscript levels.

\subsection{Operator coarsening}
Coarsening a general operator $A_l$ defined on level $l$ to level $l+1$ is performed using the well known \textit{RAT} method \cite{xu92iterative}. Specifically, to apply an operator $A_l$ to a vector $v_{l+1}$ on level $l+1$ of the mesh hierarchy: (1) the vector is interpolated onto level $l$ using the interpolation operator $L_{l+1}^l$, (2) the operator $A_l$ is applied to the interpolated vector, (3) the resulting vector is restricted back to level $l+1$ using the restriction operator $R_l^{l+1}$. This procedure is equivalent to writing a coarsened operator on mesh hierarchy level $l+1$ as
\begin{equation}
    A_{l+1} = R_l^{l+1} A_l L_{l+1}^l
\end{equation}

\section{Multigrid preconditioning}
Our plan is to utilise a multi-level $h$-multigrid solver as a right preconditioner for an iterative Krylov solver to solve the system in Equation \ref{eq:primal_form}. We use a right preconditioner instead of left since its residual is identical to the true residual. While the system matrix for Poisson's problem is symmetric positive definite, allowing for use of the conjugate gradient method, we instead opt for the GMRES algorithm as it is extendable to other problems. Furthermore, in our numerical experiments we find that convergence is generally obtained in well under 50 iterations, enabling us to consider convergence behaviour without any effects from restarts.

\subsection{Flux coarsening/Primal coarsening}
Following the discussion in \cite{saye19multigrid}, direct coarsening of the operator obtained from the primal formulation of the CDG system results in a decline in multigrid performance. Instead, each operator in the flux formulation should be individually coarsened and the Schur complement taken at each level to reform the coarse primal formulation. Coarsening of the flux formulation operator from level $l$ to level $l+1$ can be written as
\begin{equation}
    \begin{aligned}
        \begin{bmatrix} M_{l+1} & G_{l+1} \\ D_{l+1} & C_{l+1} \end{bmatrix} &=
        \begin{bmatrix} R_l^{l+1} & 0 \\ 0 & R_l^{l+1} \end{bmatrix}
        \begin{bmatrix} M_l & G_l \\ D_l & C_l \end{bmatrix}
        \begin{bmatrix} L_{l+1}^l & 0 \\ 0 & L_{l+1}^l \end{bmatrix} \\
        A_{l+1} &= C_{l+1} - D_{l+1}M^{-1}_{l+1}G_{l+1}
    \end{aligned}
\end{equation}
We verify the decline in multigrid performance from directly coarsening the primal operator in Section \ref{sect:fluxprimal}, as opposed to coarsening using the flux formulation.

\subsection{CDG switch functions}
While the CDG method has been shown to be stable and retains compactness in the primal form irrespective of the choice of switch function, it can however affect the sparsity of the matrix $C$ in Equation \ref{eq:flux_matrix_form} of the flux formulation and which can in turn affect the performance of multigrid flux operator coarsening. In particular, for each element $K_n$ separated from an element $K_i$ by a single element $K_j$, the $(K_i, K_n)$ block of the $C$ matrix is nonzero if the two conditions are satisfied:
\begin{enumerate}
    \item the switch on the edge separating $K_i,K_j$ is $S_{K_i}^{K_j} = 1$,
    \item the switch on the edge separating $K_j,K_n$ is $S_{K_j}^{K_n} = -1$.
\end{enumerate}
This implies for optimal coarsening of the operator $C$, for each partition of elements $\mathcal{T}_h$, all subsets of the partition must be closed under second neighbours that satisfy the above two properties. This is however in general impossible to satisfy for an arbitrary input mesh unless the partition consists only of one subset equal to the entire mesh.

In practice, a consistent switch function may be used to minimise the number of second neighbour interactions in $C$ not accounted for in the operator coarsening step. A consistent switch function is one where for each element $K$ with the set of neighbours $\{ K_i \}$
\begin{equation}
    | \sum_{K_i} S_K^{K_i} | < | \{ K_i \} |
\end{equation}
That is, there must be at least one inter-element boundary separating elements $K,K'$ where the switch $S_K^{K'} = -1$, and another where the switch $S_K^{K'} = 1$. The effect on performance of the multigrid preconditioner due to choice of switch function is demonstrated in Section \ref{sect:switch}.

\begin{figure}
\centering 
    \label{fig:switch_second}\includegraphics[scale=0.3]{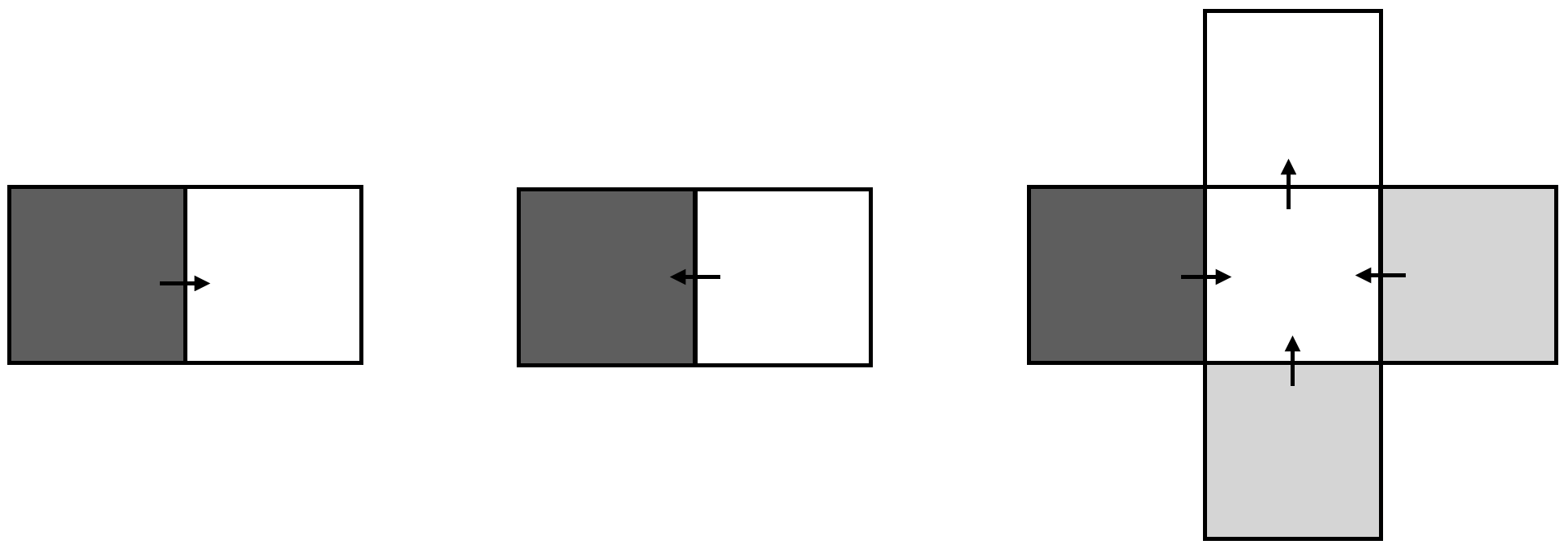}
    \caption{Effect of CDG switch function on sparsity of $C$ in Equation \ref{eq:flux_matrix_form}. On the left, arrow pointing from dark element $K_d$ to blank element $K_w$ implies switch function $S_{K_d}^{K_w} = 1$. Middle figure shows arrow pointing from blank element $K_w$ to shaded element $Kd$, denoting that $S_{K_d}^{K_w} = -1$. On the right, shaded element $K_d$ interacts with second neighbour elements shown in light grey $K_{g_1}, K_{g_2}$ as a result of the given switch, implying that blocks $(K_d, K_{g_1}), (K_d, K_{g_2})$ of the matrix C to be nonzero. }
\end{figure}

\subsection{Multigrid V-cycle}
For the h-multigrid solver, we use a single V-cycle wherein a hierarchy of meshes constructed via agglomeration is traversed using the $L^2$ projection operators outlined in Section \ref{sect:sol_transfer}. At each level, various iterations of a smoother are applied, except at the coarsest level where, the problem is solved directly.

\begin{algorithm}[H] \label{alg:vcycle}
\caption{Multigrid V-cycle}
    \begin{algorithmic}
        \STATE{Input matrix $A^{(0)}$, vector $b^{(0)}$, set $x^{(0)} = \bm{0}$}
        \STATE{Input mesh hierarchy $E = \{ E_0, E_1,...,E_n \}$}
        \FOR{$k = 1:n-1$}
            \STATE{Project to current level $b^{k} = R_{k-1}^{k} ( b^{k-1} - A^{k-1}x^{k-1})$}
            \STATE{Construct smoother $\tilde{A}_{k}$ from $A_{k}$}
            \FOR{$i = n_{pre}$}
                \STATE{Apply smoother, $x^k = x^k + \alpha \tilde{A}_{k}^{-1} (b^{k} - A^{k}x^k)$}
            \ENDFOR
        \ENDFOR
        \STATE{Solve $A^n x^n = b^n$ directly}
        \FOR{$k = n-1:1$}
            \STATE{Project to current level $x^{k} = x^{k} + L_{k+1}^{k} x^{k+1} $}
            \STATE{Construct smoother $\tilde{A}_{k}$ from $A_{k}$}
            \FOR{$i = n_{pre}$}
                \STATE{Apply smoother, $x^k = x^k + \alpha \tilde{A}_{k}^{-1} (b^{k} - A^{k}x^k)$}
            \ENDFOR
        \ENDFOR
    \RETURN $x^{(0)}$
    \end{algorithmic}
\end{algorithm}
Commonly used smoother include block Jacobi, block Gauss-Seidel, or incomplete LU factorisations. In this work we focus on block Jacobi smoothers with a damping factor $\alpha = \frac{2}{3}$, as they are simple to parallelise for large systems.

\subsection{$hp$-multigrid}
\label{sect:hp}
In this manuscript, we consider explicitly only the case of linear basis functions in the Discontinuous Galerkin discretisation. For problems with higher polynomial degree basis functions, we would
first employ standard $p$-multigrid \cite{fidkowski05multigrid,persson08newtongmres} on the
fine mesh to project down to $p=1$ basis functions, and then followed by our agglomeration method.
For simplicity and to only highlight the $h$-multigrid procedure, we only consider $p=1$
in our examples.

\section{Numerical results}
In this section we present numerical results to evaluate the performance of the multigrid preconditioner. Unless otherwise stated, a consistent switch function is used for flux definition in the CDG discretization. Our initial solution vector is always set as the zero vector, and we iterate until we reach a tolerance of $10^{-8}$ in the relative norm $|Au_h - b| / |b|$. In all examples, we agglomerate elements in the mesh until the lowest level in the $h$-multigrid hierarchy consists of only one element using the greedy algorithm described above.

Following the discussion in Section \ref{sect:hp}, for all the following examples, we consider only linear order basis functions for both basis functions defined on the initial mesh and on all subsequent meshes in the $h$-multigrid hierarchy. Within the multigrid V-cycle we choose the number of pre-smoothing steps $n_{pre}=0$, and the number of post-smoothing steps $n_{post}=3$. We choose a GMRES restart parameter of 50, as in all the following examples we converge in fewer iterations and so convergence is not impacted by any restarts. We also do not consider true computational time in this study, and report only the number of iterations required for convergence.

\subsection{Flux vs primal coarsening and choice of Dirichlet parameter}
\label{sect:fluxprimal}
We start by solving Poisson's problem on the domain $\Omega = [0,1]^2$ using a uniform square $n \times n$ mesh. We impose Neumann conditions on the vertical boundaries at $x=0,1$, in addition to Dirichlet boundary conditions on the horizontal boundaries at $y=0,1$. We build the h-multigrid hierarchy using Algorithm \ref{alg:agglomeration}, which is shown in the top row of Figure \ref{fig:square_hgrid}.

\begin{figure}
    \centering
    \subfloat[Flux coarsening]{%
        \resizebox{0.4\linewidth}{!}{%
            \begin{tikzpicture}[font=\Large]
                \begin{semilogxaxis}[
                	xlabel = {$n$},
                	ylabel = {GMRES Iterations},
                	ymin=0, ymax=30,
                	xmin=1, xmax=64,
                    xtick={2,4,8,16,32},
                    xticklabels={2,4,8,16,32},
            		grid = both,
                	grid style = {line width=.1pt, draw=gray!15},
            		major grid style = {line width=.2pt, draw=gray!50},
                ]
                \addplot[mark=o, color=black, mark size = 3pt]
                table[ x=h, y=f1 ]{dat/quadsquare_iters.csv};
                
                \addplot[mark=x, color=black, mark size =3pt]
                table[ x=h, y=f2 ]{dat/quadsquare_iters.csv};
                
                \addplot[mark=*, color=black, mark size=3pt]
                table[ x=h, y=f3 ]{dat/quadsquare_iters.csv};
                
                \addplot[mark=pentagon, color=black, mark size=3pt]
                table[ x=h, y=f4 ]{dat/quadsquare_iters.csv};
                \end{semilogxaxis}
            \end{tikzpicture}%
        }
    }\hfil
    \subfloat[Primal coarsening]{%
        \resizebox{0.4\linewidth}{!}{%
            \begin{tikzpicture}[font=\Large]
                \begin{semilogxaxis}[
                	xlabel = {$n$},
                	ylabel = {GMRES Iterations},
                	ymin=0, ymax=30,
                	xmin=1, xmax=64,
                    xtick={2,4,8,16,32},
                    xticklabels={2,4,8,16,32},
            		grid = both,
                	grid style = {line width=.1pt, draw=gray!15},
            		major grid style = {line width=.2pt, draw=gray!50},
                ]
                \addplot[mark=o, color=black, mark size = 3pt]
                table[ x=h, y=p1 ]{dat/quadsquare_iters.csv};
                \label{plot:quadsquare1}
                
                \addplot[mark=x, color=black, mark size = 3pt]
                table[ x=h, y=p2 ]{dat/quadsquare_iters.csv};
                \label{plot:quadsquare2}
                
                \addplot[mark=*, color=black, mark size = 3pt]
                table[ x=h, y=p3 ]{dat/quadsquare_iters.csv};
                \label{plot:quadsquare3}
                
                \addplot[mark=pentagon, color=black, mark size = 3pt]
                table[ x=h, y=p4 ]{dat/quadsquare_iters.csv};
                \label{plot:quadsquare4}
                \end{semilogxaxis}
            \end{tikzpicture}%
        }
    }
    
\caption{h-multigrid convergence for flux vs. primal coarsening on a square $n \times n$ mesh. Different plot markers indicate varying Dirichlet parameters $C_D$: \ref{plot:quadsquare1}, \ref{plot:quadsquare2}, \ref{plot:quadsquare3}, \ref{plot:quadsquare4} denote values of $C_D = 10^1/h_{avg},~ 10^2/h_{avg},~ 10^3/h_{avg},~ 10^4/h_{avg},$ respectively. }
\label{fig:quadsquare_iter}
\end{figure}
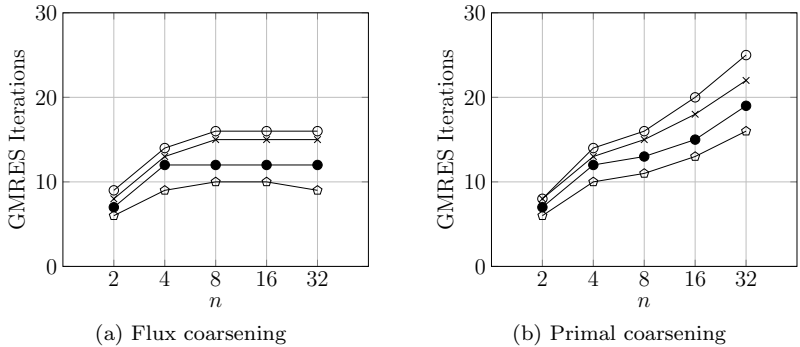

Figure \ref{fig:quadsquare_iter} shows the number of iterations to convergence for the square mesh under h-refinement, in addition to varying values of the Dirichlet penalty parameter $C_D$. A deterioration in performance under primal coarsening as $h_{avg} \rightarrow 0$ is seen as in \cite{saye19multigrid}, whereas performance under flux coarsening does not suffer similar problems 

Increased performance is also gained through using larger values of $C_D$, a result the choice of block Jacobi as the smoother in the h-multigrid solver. As only the magnitude of values in the blocks on the diagonal of matrix $C$ in Equation \ref{eq:flux_matrix_form} scale with the value of $C_D$, an increase in $C_D$ implies an increase to the values in the blocks on the diagonal of $A$ in Equation \ref{eq:primal_form} relative to values in blocks off the diagonal of $A$. Based on these observations, going forward for the remainder of our tests, we focus on flux coarsening using a value of $C_D = 10^4/h_{avg}$.

\subsubsection{Hierarchy element shapes}
\begin{figure}
    \centering
    \includegraphics[scale=0.6]{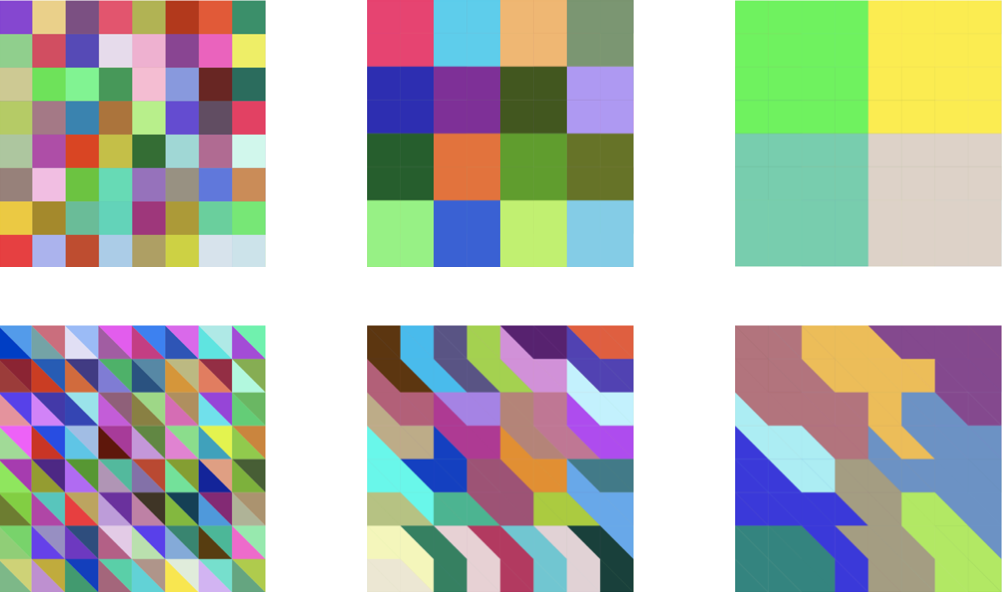}
    \caption{Regular and irregular mesh hierarchies. Mesh hierarchies shown from left to right correspond to level = 0,1,2 of a $8 x 8$ square mesh respectively.}
    \label{fig:square_hgrid}
\end{figure}

We investigate the effect of irregular element shapes in the h-multigrid hierarchy on performance by considering once more a domain $\Omega = [0,1]^2$, with Neumann conditions on the boundaries at $x=0,1$, and Dirichlet conditions on the boundaries at $y=0,1$. We compare the number of iterations to convergence using two different mesh hierarchies as shown in Figure \ref{fig:square_hgrid}, one with regular quadrilateral shaped elements at each level in the h-multigrid and the other with highly irregularly shaped elements that are in general non-convex.

\begin{figure}
    \centering
    \resizebox{0.5\linewidth}{!}{%
        \begin{tikzpicture}[font=\Large]
            \begin{semilogxaxis}[
                xlabel = {$n$},
            	ylabel = {GMRES Iterations},
            	ymin=0, ymax=30,
            	xmin=1, xmax=64,
                xtick={2,4,8,16,32},
                xticklabels={2,4,8,16,32},
            	grid = both,
            	grid style = {line width=.1pt, draw=gray!15},
        		major grid style = {line width=.2pt, draw=gray!50},
            ]
            \addplot[mark=square, color=black, mark size=3pt]
            table[ x=h, y=q ]{dat/triquadsquare_iters.csv};
            \label{plot:regular_square}
                
            \addplot[mark=triangle, color=black, mark size = 3pt]
            table[ x=h, y=t ]{dat/triquadsquare_iters.csv};
            \label{plot:irregular_square}
            \end{semilogxaxis}
        \end{tikzpicture}%
        \label{fig:square_comparison}
    }
\caption{Comparison of multigrid performance with varying mesh hierarchy element shapes on a square $n \times n$ mesh. \ref{plot:regular_square} denotes the number of iterations on a regular mesh hierarchy, \ref{plot:irregular_square} denotes the number of iterations on a non-regular mesh hierarchy.}
\label{fig:triquad_iter}
\end{figure}
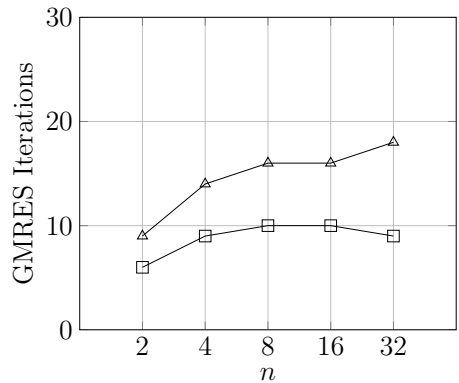

The plot in Figure \ref{fig:square_comparison} show the number of iterations to convergence using the regular and irregular shaped elements respectively. While a decrease in performance is observed in using irregularly shaped elements in the h-multigrid hierarchy, the decrease is independent of element size $h_{avg}$.

\subsection{NACA airfoil} \label{sect:NACA}
\begin{figure}
    \centering 
    \includegraphics[scale=0.35]{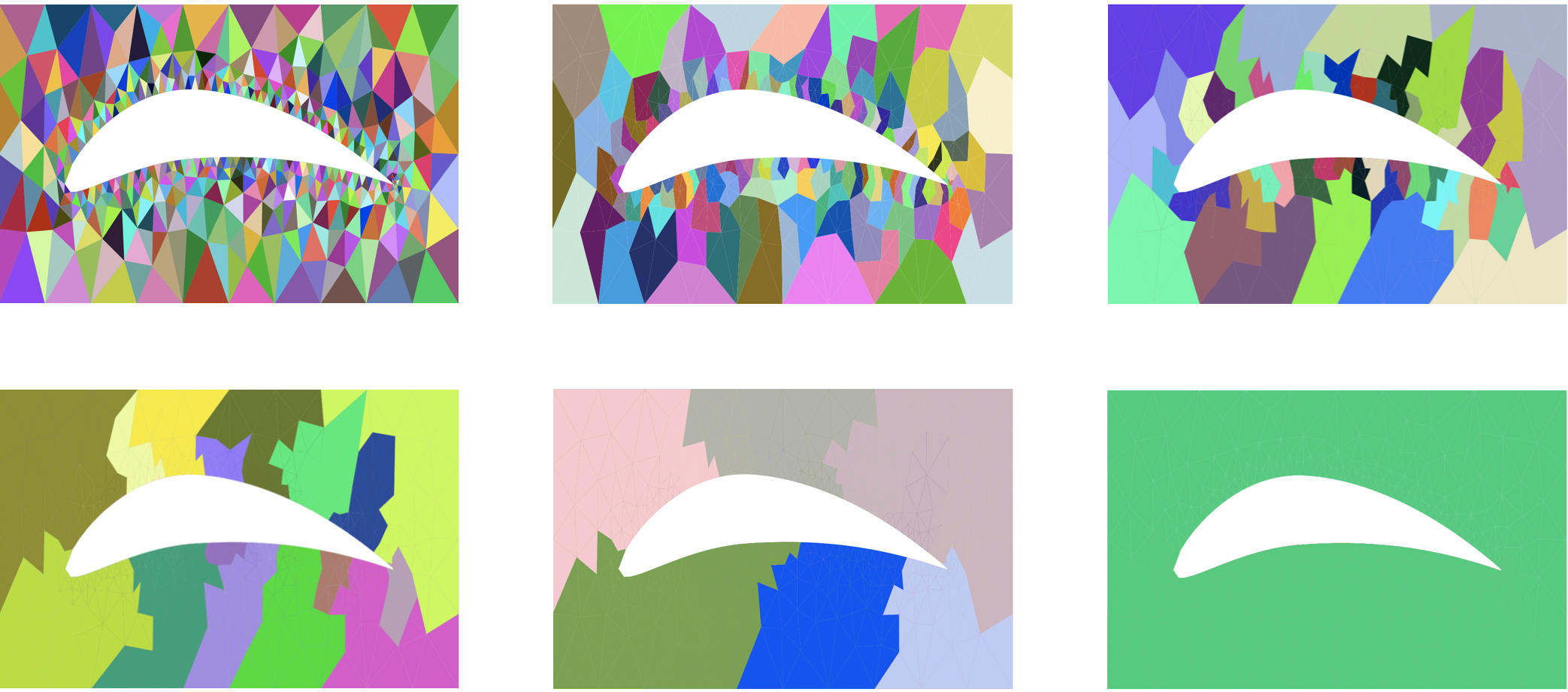}
    \caption{Multigrid hierarchy for airfoil mesh.}
\label{fig:airfoil_h}
\end{figure}

To investigate the effect of non-uniform element sizes and also of meshes which are not simply connected, we consider the example of Poisson's problem on a rectangular domain around a NACA airfoil. Figure \ref{fig:airfoil_h} shows the h-multigrid hierarchy of the coarsest input mesh consisting of 605 elements, shown in the top left of the figure. Dirichlet conditions are applied on the boundary at the airfoil and at the two horizontal boundaries, while Neumann conditions are applied at the two vertical boundaries.

\subsubsection{Switch function} \label{sect:switch}
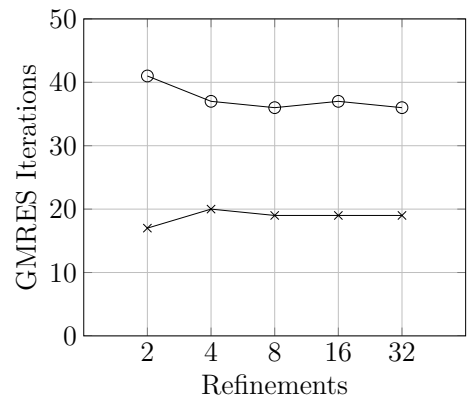
\begin{figure}
    \centering
    \resizebox{0.5\linewidth}{!}{%
    \begin{tikzpicture}[font=\Large]
        \begin{semilogxaxis}[
            xlabel = {Refinements},
            ylabel = {GMRES Iterations},
            ymin=0, ymax=50,
            xmin=1, xmax=64,
            xtick={2,4,8,16,32},
            xticklabels={2,4,8,16,32},
            grid = both,
            grid style = {line width=.1pt, draw=gray!15},
            major grid style = {line width=.2pt, draw=gray!50},
        ]
        \addplot[mark=x, color=black, mark size=3pt]
        table[ x=h, y=Consistent ]{dat/airfoil_iters.csv};
        \label{plot:airfoil_consistent}
                
        \addplot[mark=o, color=black, mark size=3pt]
        table[ x=h, y=Natural ]{dat/airfoil_iters.csv};
        \label{plot:airfoil_natural}
        \end{semilogxaxis}
    \end{tikzpicture}%
}
\caption{GMRES convergence for Poisson's problem on airfoil mesh. \ref{plot:airfoil_consistent} shows the number of iterations using a consistent switch function, \ref{plot:airfoil_natural} shows the number of iterations using a natural switch function.}
\label{fig:airfoil_iters}
\end{figure}

Figure \ref{fig:airfoil_iters} shows the number of iterations to convergence for the airfoil problem using a consistent switch function, and a natural switch function based on random element enumeration. The performance of the multigrid preconditioner is shown to clearly deteriorate with a poorly chosen switch function. Using a consistent switch function, the presence of non simply connected elements in the mesh hierarchy does not seem to have a large effect on the performance of the multigrid preconditioner.

\subsection{Convection-Diffusion} \label{sect:switch}

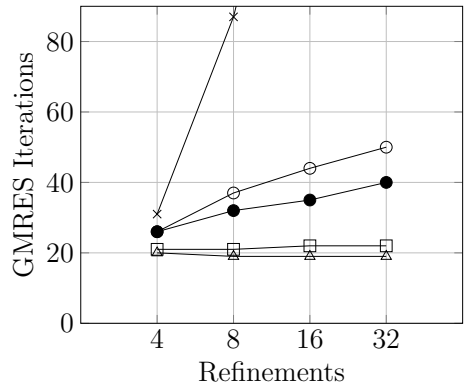
\begin{figure}
    \centering
    \resizebox{0.5\linewidth}{!}{%
    \begin{tikzpicture}[font=\Large]
        \begin{semilogxaxis}[
            xlabel = {Refinements},
            ylabel = {GMRES Iterations},
            ymin=0, ymax=90,
            xmin=2, xmax=64,
            xtick={4,8,16,32},
            xticklabels={4,8,16,32},
            grid = both,
            grid style = {line width=.1pt, draw=gray!15},
            major grid style = {line width=.2pt, draw=gray!50},
        ]
        \addplot[mark=x, color=black, mark size=3pt]
        table[ x=h, y=c100 ]{dat/airfoil_convdiff.csv};
        \label{plot:airfoil_cd_100}
                
        \addplot[mark=o, color=black, mark size=3pt]
        table[ x=h, y=c10 ]{dat/airfoil_convdiff.csv};
        \label{plot:airfoil_cd_10}
        
        \addplot[mark=*, color=black, mark size=3pt]
        table[ x=h, y=c1 ]{dat/airfoil_convdiff.csv};
        \label{plot:airfoil_cd_1}
        
        \addplot[mark=square, color=black, mark size=3pt]
        table[ x=h, y=c01 ]{dat/airfoil_convdiff.csv};
        \label{plot:airfoil_cd_01}
        
        \addplot[mark=triangle, color=black, mark size=3pt]
        table[ x=h, y=c0 ]{dat/airfoil_convdiff.csv};
        \label{plot:airfoil_cd_0}
        
        \end{semilogxaxis}
    \end{tikzpicture}%
}
\caption{GMRES convergence for convection-diffusion equation on airfoil mesh. The number of iterations are shown for $\beta = 100$ by \ref{plot:airfoil_cd_100}, $\beta = 10$ by \ref{plot:airfoil_cd_10}, $\beta = 1$ by \ref{plot:airfoil_cd_1}, $\beta = 0.1$ by \ref{plot:airfoil_cd_01} and $\beta = 0$ by \ref{plot:airfoil_cd_0}. }
\label{fig:airfoil_convdiff}
\end{figure}

Finally, we consider the more general example of a convection-diffusion equation on
the airfoil mesh:
\begin{equation}
    \beta \bm{v} \cdot \nabla u + \Delta u = f.
\end{equation}
We employ zero Dirichlet boundary conditions everywhere, we set $f=1$, and the velocity
field $\bm{v}=(1,0)$. The resulting convergence in the GMRES iterations is shown in Figure~\ref{fig:airfoil_convdiff}, for a range of values of $\beta$ under refinement.

We see that the performance of the preconditioner is largely unaffected for small values of $\beta$, but quickly deteriorates with higher magnitudes of $\beta$. This is expected as it changes the structure of the problem. We note however that this can be
fixed by using other existing smoothers for convection such as line-based
solvers \cite{fidkowski05multigrid}, or ILU/Gauss-Seidel with good element ordering \cite{persson08newtongmres}.

\section{Conclusions}
\label{sec:conclusions}

We have developed an algorithm for constructing suitable mesh hierarchies for the geometric multigrid method via use of simple element agglomeration. The merged elements will in general be polyhedral, which are easily supported using a discontinuous Galerkin discretization. While the method should perform well with any choice of numerical fluxes, we have used the Compact DG method and showed that in this case a consistent switch function gives better multigrid performance. The resulting solver gives excellent performance for Poisson's equation on fully unstructured meshes, as well as for convection-diffusion with moderate magnitudes of the convective component. Future work include extension to other equations, parallelization, and numerical examples in 3D.

\section*{Acknowledgments}

This work was supported in part by the Director, Office of Science, Office of
Advanced Scientific Computing Research, U.S. Department of Energy under
Contract No. DE-AC02-05CH11231.

\bibliographystyle{plain}
\bibliography{references}

\end{document}